\documentstyle[12pt]{article}
\oddsidemargin=\evensidemargin
\newfont{\blb}{msbm10 scaled\magstep1}

\newtheorem{theo}{Theorem}[section]

\newtheorem{prop}[theo]{Proposition}
\newtheorem{lemm}[theo]{Lemma}

\begin{document}

\date{}
\author{Anastasia Hadjievangelou \& Gunnar Traustason \\
Department of  Mathematical Sciences, \\
University of Bath, \\
Bath BA2 7AY,
UK}

\title{Sandwich groups and (strong) left $3$-Engel elements in groups}

\maketitle
\begin{abstract}
\mbox{}\\
In this paper we prove a group theoretic analogue of the well known local nilpotence theorem for sandwich Lie algebras due to Kostrikin and Zel'manov. We introduce the notion of a strong left $3$-Engel element of a group $G$ and show that these are always in the locally nilpotent radical of $G$. This generalises a previous result of Jabara and Traustason that showed that a left $3$-Engel element $a$ of a group $G$ is in the locally nilpotent radical of $G$ whenever $a$ is of odd order.\footnote{Primary 20F45, Seondary 20F12}
\end{abstract}

\section{Introduction}
\subsection{Some background}

Let $G$ be a group. An element $a\in G$ is a left Engel element in $G$, if for
each $x\in G$ there exists a non-negative integer $n(x)$ such that
      $$[[[x,\underbrace {a],a],\ldots ,a]}_{n(x)}=1.$$
If $n(x)$ is bounded above by $n$ then we say that $a$ is a left $n$-Engel element
in $G$. It is straightforward to see that any element of the Hirsch-Plotkin
radical $HP(G)$ of $G$ is a left Engel element and the converse is known
to be true for some classes of groups, including solvable groups and 
finite groups (more generally groups satisfying the maximal condition on
subgroups) \cite{Gru, Baer}. The converse is however not true in general and this is the case
even for bounded left Engel elements. In fact whereas one sees readily that 
a left $2$-Engel element is always in the Hirsch-Plotkin radical this
is still an open question for left  $3$-Engel elements in general. There has though been some breakthrough in recent years and in \cite{Jab} it is shown that any left $3$-Engel element of odd order is contained in $HP(G)$. In this paper we will generalise this result to include elements of any order, but we need to replace left $3$-Engel with a stronger condition that we call strong left $3$-Engel, that is however equivalent to left $3$-Engel when the element is of odd order. \\ \\
From \cite{Trac} one knows that in order to show that any left $3$-Engel element of finite order is in the Hirsch-Plotkin radical, it suffices to deal with elements order $2$. From looking at similar setting for Lie algebra (see section 1.2), there are reasons to doubt that this is true however. Much less is known about left $4$-Engel elements, although there are some interesting results in \cite{Ab2}.  \\ \\
Groups of prime power exponent are known to satisfy some Engel type conditions and the solution to the restricted Burnside problem in particular makes use of the fact that the associated Lie ring satisfies certain Engel type identities \cite{zc,zd}. Considering left Engel elements, it was observed by William Burnside \cite{Burn} that every element in a group of exponent 3 is a left 2-Engel element and so the fact that every left 2-Engel element lies in the Hirsch-Plotkin radical can be seen as the underlying reason why groups of exponent 3 are locally finite. For groups of 2-power exponent there is a close link with left Engel elements. Let  $G$ be a finitely generated group of exponent $2^{n}$ and $a$ an element in $G$ of order $2$, then 
            $$[[[x,\underbrace{a],a],\ldots ,a}_{n+1}]=[x,a]^{(-2)^{n}}=1.$$ 
Thus $a$ is a left $(n+1)$-Engel element of $G$. It follows from this that
if $G/G^{2^{n-1}}$ is finite and the left $(n+1)$-Engel elements of $G$ are
in the Hirsch-Plotkin radical, then $G$ is finite. As we know that for sufficiently large $n$ the variety of groups of exponent $2^{n}$ is not locally finite \cite{Ivan,Lys}, it follows that for sufficiently large $n$ there are left $n$-Engel elements that are not contained in the Hirsch-Plotkin radical. Using the fact that groups of exponent 4 are locally finite \cite{San}, one can also see that if all left 4-Engel elements of a group $G$ of exponent 8 are in $HP(G)$ then $G$ is locally finite. \\ \\
Swapping the role of $a$ and $x$ in the definition of a left Engel element we get the notion of a right Engel element. Thus an element $a\in G$ is a right Engel element, if for each $x\in G$ there exists a non-negative 
integer $n(x)$ such that 
    $$[a,_{n(x)} x]=1.$$
If $n(x)$ is bounded above by $n$, we say that $a$ is a right $n$-Engel element. By a classical result
of Heineken \cite{Hein1} one knows that if $a$ is a right $n$-Engel element in $G$ then $a^{-1}$ is a left $(n+1)$-Engel
element. \\ \\
In \cite{New} M. Newell proved that if $a$ is a right $3$-Engel element in $G$ then $a\in HP(G)$ and in
fact he proved the stronger result that $\langle a\rangle^{G}$ is nilpotent of class at most $3$. The natural question arises whether the analogous result holds for left $3$-Engel elements. In \cite{GGM} it has been shown that this is not the case by giving an example of a locally finite $2$-group with a left $3$-Engel element $a$ such that $\langle a\rangle^{G}$ is not nilpotent and in \cite{GA} this was generalised to an infinite family of examples. In \cite{GAM} this was extended by giving an example, for any odd prime $p$, of a locally finite $p$-group containing a left $3$-Engel element $a$ such that $\langle a\rangle^{G}$ is not nilpotent.
\subsection{(Strong) Left $3$-Engel elements and sandwich groups}
The approach here for (strong) left $3$-Engel elements, as in \cite{Jab},  is through working with  sandwich groups, but these can be seen as the group theoretic analogue of  sandwich Lie algebras that were introduced by Kostrikin \cite{Ko1}. Kostrikin and Zel'manov proved that sandwich Lie algebras are locally nilpotent \cite{Kos} and this fact is a key ingredient to both Kostrikin's solution to the restricted Burnside problem for groups of prime exponent \cite{Ko2}, and Zel'manov's general solution \cite{zc,zd}. Let us recall the definition.  As for group commutators we are using left normed notion for Lie products. \\ \\
{\bf Definition}. Let $L$ be a Lie algebra and $a \in G$. We say that $a$ is a {\it sandwich element} of $L$,  if  $axa=0$ and $axya=0$ for all $x,y \in L$. \\ \\
{\bf Remark}.  It is not difficult to see that the latter condition is superfluous when the characteristic is odd.  This follows from $0=x(yaa)=xyaa-2xaya+xaay=2axya$. \\ \\
{\bf Definition}. A Lie algebra $L$ is a sandwich algebra if it can be generated (as Lie algebra) by sandwich elements. \\ \\
The reason why the latter condition, in the definition of a sandwich element, is needed can be seen from the following non-nilpotent  elementary example of a $3$-generator Lie algebra over $\mbox{GF}(2)$, where the three generators only satisfy the first condition of a sandwich element. \\ \\ 
{\bf Example}.  Consider the largest Lie algebra $L= \langle a,b,c \rangle$ over $GF(2)$ subject to $Id(c)$ being abelian, $bc=0$ and $bxb=axa=cxc=0$, for all $x \in L$. One readily sees that $L$ is non-nilpotent with basis $a$, $b$, $ab$, $c(ab)^{n},\ n\geq 0$, $c(ab)^{n}a,\ n\geq 0$. \\ \\
In \cite{Trau}  the notion of a sandwich group was introduced by the second author. \\ \\
{\bf Definition}. Let $G$ be a group and $X$ a subset of $G$. We say that $X$ is a {\it  sandwich set},  if  $\langle a,b^{g}\rangle$ is nilpotent of class most $2$ for all $a,b\in X$
and $g\in \langle X \rangle$. We say that a group $G$ is a {\it sandwich group} if $G=\langle X\rangle$, where $X$ is a sandwich set in $G$.  \\ \\
The connection with left $3$-Engel elements comes from the fact that the following are equivalent.   \vspace{2 mm} \\
(1) For every pair $(G, a)$ where $a$ is a left 3-Engel element in the group
$G$ we have that $a$ is in the locally nilpotent radical of $G$.\\
(2) Every finitely generated sandwich group is  nilpotent.\\ \\
The reason for this is that firstly $\langle a\rangle^{G}$ is a sandwich group whenever $a$ is a left $3$-Engel element in $G$, and secondly any element $a$ of a sandwich set $X$ is left $3$-Engel in $\langle X\rangle$. The advantage of working in sandwich groups is that there is a largest one for each given rank and so we have something specific to work with that does not rely on sitting inside a larger group. \\ \\
In order to show that a left $3$-Engel element $a$ in $G$ is in the Hirsch-Plotkin radical of $G$, when $a$ is of odd order, it was shown in \cite{Jab} that any sandwich group, generated a finite sandwich set consisting of elements of odd order, is nilpotent. \\ \\ 
We would like to extend this result to include groups generated by elements of any order. Here we can take our cue from the definition of sandwich algebras. The result is the following. \\ \\
{\bf Definition}. Let $G$ be a group and $X$ a subset of $G$. We say that $X$ is a {\it strong sandwich set},  if  \vspace{2 mm} \\
(1) $\langle a,b^{g}\rangle$ is nilpotent of class most $2$ for all $a,b\in X$
and $g\in \langle X \rangle$.\\
(2) $\langle a,b^{f},c^{g}\rangle$ is nilpotent of class at most $3$ for all $a,b,c\in X$ and $f,g\in \langle X\rangle$. \vspace{2mm} \\
We say that a group $G$ is a strong sandwich group, if it can be generated by a strong sandwich set. \\ \\
{\bf Remark}. We have here something analogous to sandwich Lie algebras. Namely that if all the elements in the generating set $X$ are of odd order then the 2nd condition is superfluous. This follows from the fact that a sandwich group of rank $3$ is nilpotent of class at most $3$ if all the elements are of odd order \cite{Trau}. The free sandwich group of rank $3$ is however nilpotent of class $5$. \\ \\
Here there is also a connection with left $3$-Engel elements. We first need a definition.  \\ \\
{\bf Definition}. We say that an element $a\in G$ is a strong left $3$-Engel element in $G$  if  \vspace{2 mm} \\
(1) $\langle a,a^{g}\rangle$ is nilpotent of class at most $2$ for all  $g\in G$.\\
(2) $\langle a,a^{f},a^{g}\rangle$ is nilpotent of class at most $3$ for all $f,g\in G$.\\ \\
{\bf Remark}. Notice that $a$ is left $3$-Engel in $G$ if and only if it satisfies (1). The condition (2) is superfluous when $a$ is of odd order and therefore a left $3$-Engel element of odd order is a strong left $3$-Engel element. \\ \\
Strong left $3$-Engel elements and strong sandwich groups are related through the fact that the following are equivalent: \vspace{2 mm} \\
(1) $\langle a \rangle ^G$ is locally nilpotent whenever $a$ is a strong left $3$-Engel element in $G$. \\
(2) Every finitely generated strong sandwich group  is nilpotent.\\ \\
The main theorems of this paper are the following generalisations of the results of Jabara and Traustason \cite{Jab}.  
\begin{theo}
Every finitely generated strong sandwich group is nilpotent.
\end{theo}
\begin{theo}
If $a$ is strong left $3$-Engel element of a group $G$, then $\langle a\rangle^{G}$ is locally nilpotent.
\end{theo}
{\bf Remark}. Whereas we have seen that there is an elementary example that shows that the 2nd condition in the definition of a sandwich Lie algebra is needed, it is still an open question whether a left $3$-Engel element in a group $G$ must always be in the Hirsch-Plotkin radical. This is even the case under the extra condition that the group is residually finite. Of course we know this is true if the element is of odd order. \\ \\
We end this introduction by mentioning two applications of the Theorem of Jabara and Traustason for groups of exponents $5$ and $9$. \\ \\
{\bf Theorem A} \cite{Trau}. {\it Let $G$ be a group of exponent $5$. Then $G$ is locally finite if and only if it satisfies the law} 
     $$[z,[y,x,x,x],[y,x,x,x],[y,x,x,x]]=1.$$
{\bf Remark}. This implies in particular that a group of exponent $5$ is locally finite if and only if all the $3$-generator subgroups are finite. This was originally proved by Vaughan-Lee \cite{Vaug1}. \\ \\
{\bf Theorem B} \cite{Jab}. {\it Let $w$ be a law in $n$ variables $x_{1},\ldots ,x_{n}$ where the variety of groups satisfying the law $w^{3}=1$ is a locally finite variety of groups of exponent $9$. Then the same is true for the variety of groups satisfying the law $(x_{n+1}^{3}w^{3})^{3}=1$.} \\ \\
{\bf Remark}. We can use Theorem B to come up with an explicit sequence of words. Define the word $W_{n}=W_{n}(x_{1},\ldots , x_{n})$ in $n$ variables recursively by $w_{1}=x_{1}$
and $w_{n+1}=x_{n+1}^{3}w_{n}^{3}$. The variety of groups satisfying the law $x_{1}^{3}=1$ is locally finite by Burnside and by repeated application of Theorem B we see that, for each $n\geq 1$, the variety of groups satisfying the law $w_{n}^{3}$ is a locally finite variety of groups of exponent $9$.

\section{Commutator closure of strong sandwich sets}
We want to show that strong sandwich sets are closed under taking commutators. In order to show this we first deal with a weaker version of sandwich sets. \\ \\
{\bf Definition}. Let $G$ be a group and $X$ a subset of $G$. We say that $X$ is a {\it partial strong sandwich set},  if  \vspace{2 mm} \\
(1) $\langle a,b^{g}\rangle$ is nilpotent of class most $2$ for all $a,b\in X$, where $a$ and $b$ are distinct,
and $g\in \langle X \rangle$.\\
(2) $\langle a,b^{f},c^{g}\rangle$ is nilpotent of class at most $3$ for all $a,b,c\in X$, where $a,b$ and $c$ are distinct,  and $f,g\in \langle X\rangle$. 
% 
%
%Let $F$ be the free group of rank $r$ generated freely by $x_{1},\ldots ,x_{r}$. Let $R$ be the smallest normal subgroup of $F$ for which $\{ x_{1}R,\ldots, x_{r}R\}$ be comes a strong sandwich set. We call $F/R$ the {\it free strong sandwich group of rank $r$}.
%
\begin{prop}\label{prop1}
Let $G=\{a,b,c,d\}$ be a partial strong sandwich set of some group. Then $G$ is nilpotent of class at most $5$. 
\end{prop}
In order to prove this proposition we will need the following lemma.
\begin{lemm}\label{lem1}
We have that H=$\langle a,a^{b},c,d\rangle$ is nilpotent of class at most $4$.
\end{lemm}
{\bf Proof}\ \ We prove this in few steps. Notice that we can only use the weak partial sandwich properties for (1) and (2).  \\ \\
\underline{Step 1}. We have that  $[a^{b},c,[d,a]]$ is in $Z(H)$. \\ \\
To see this, notice first that $[a^{b},c]$ commutes with $a^{b}$ by (1) and, as $[a^{b},c]=[a[a,b],c]=[a,c]^{[a,b]}[a,b,c]$, it also commutes with $a$ by (1) and (2). Also $[d,a]$ commutes with $a$ by (1) and with $a^{b}=a[a,b]$  by (1) and (2). This shows that $[a^{b},c,[d,a]]$ commutes with both $a$ and $a^{b}$. \\ \\
Then
 \begin{eqnarray*}
\mbox{}[a^{b},c,[d,a]] & = & [a^{b},c,a^{-d}a] \\
                                       & = & [a^{b},c,aa^{-d}] \mbox{\ \ (as }a^{d}=a[a,d]\mbox{ commutes with }a\mbox{ by (1))} \\
                                      & = & [a^{b},c,a^{-d}] \mbox{\ \ (as }[a^{b},c] \mbox{ commutes with }a) 
\end{eqnarray*}  
But, as $[a^{-1}a^{b},c,a^{-d}]=[[a^{-1},c]^{a^{b}}[a^{b},c],a^{-d}]=[a^{b},c,a^{-d}]$ (using $[a^{-1},c]$ commutes with $a^{b}=a[a,b]$  and $a^{d}=a[a,d]$ by (1) and (2)), we see from this and the above that 
 \begin{eqnarray*}
\mbox{}[a^{b},c,[d,a]] & = & [a,b,c,a^{-d}] \\
                                      & = & [b^{-a}b,c,a^{-d}] \\ 
                                      & = & [bb^{-a},c,a^{-d}] \mbox{\ \ (as }b^{a}=b[b,a]\mbox{ commutes with }b\mbox{ by (1))} \\
                                     & = & [[b,c][b^{-a},c],a^{-d}] \mbox{\ \ (as }[b,c]\mbox{ commutes with }b^{a}=b[b,a]\mbox{ by (1,2).}  \\
                                   & = & [b,c,a^{-d}]^{[b^{-a},c]}[b^{-a},c,a^{-d}].                                 
\end{eqnarray*}
This last element commutes with $c$ by (1) and (2) and thus we have seen that $[a^{b},c,[d,a]]$ commutes with $c$. In order to show that $[a^{b},c,[d,a]]$ is in $Z(H)$, it remains to see that
it commutes with $d$. As  
      $$[a^{b},c,[d,a]]=[c^{-a^{b}}c,[d,a]]=[c^{-a^{b}},[d,a]][c^{-a^{b}},[d,a],c][c,[d,a]].$$         
It follows from (2) that it suffices to show that $[c^{-a^{b}},[d,a],c]$ commutes with $d$. Using (1) and (2) again, we know that $c^{-a^{b}}=c^{-1}[c^{-1},a^{b}]$, $[d,a]$ and $c$ all commute with $[c,d]$ and $[c^{-a^{b}},[d,a]]$ commutes with $d$. Hence 
\begin{eqnarray*}
     [c^{-a^{b}},[d,a],c]^{d} & = & [c^{-a^{b}},[d,a],c]^{d} = [c^{-a^{b}},[d,a],c[c,d]] \\
                              & = & [c^{-a^{b}},[d,a],c]^{[c,d]} = [c^{-a^{b}},[d,a],c]
\end{eqnarray*}
and we have shown that $[c^{-a^{b}},[d,a],c]$ commutes with $d$ and hence $[a^{b},c,[d,a]]$ is in $Z(H)$. \\ \\
\underline{Step 2}. $[a^{b},c,d]$ is in $Z_{2}(H)$.  \\ \\ 
Notice that $d$ commutes with $[d,a]$ by (1) and that $[a^{b},c]=[a[a,b],c]=[a,c]^{[a,b]}[a,b,c]$ commutes with $a$ by (1) and (2).  Using this and Step 1, we see that modulo $Z(H)$
         $$[a^{b},c,d]^{a}=[a^{b},c,d[d,a]]=[a^{b},c,d]^{[d,a]}=[a^{b},c,d].$$
Thus $[a^{b},c,d]$ commutes with $a$ modulo $Z(H)$. By (2) it also commutes with $a^{b},c$ and $d$. Hence $[a^{b},c,d]$ is in $Z_{2}(H)$. \\ \\
\underline{Step 3}. $[a^{b},c]$ and $[a^{b},d]$ are in $Z_{3}(H)$. \\ \\ 
We know from Step 2 that $[a^{b},c]$ commutes with $d$ modulo $Z_{2}(H)$. By (1) it also commutes with $a^{b}$ and $c$. Finally $[a^{b},c]=[a[a,b],c]=[a,c]^{[a,b]}[a,c]$ commutes with
$a$ by (1) and thus is in $Z_{3}(H)$. By symmetry $[a^{b},d]$ is also in $Z_{3}(H)$. \\ \\
\underline{Step 4}. $H$ is nilpotent of class at most $4$.  \\ \\
As $a^{b}=a[a,b]$ commutes with $a$ by (1) and obviously with $a^{b}$, it follows from Step 3 that $a^{b}$ is in $Z_{4}(G)$. Hence $\gamma_{5}(H)=\gamma_{5}(\langle a,b,c\rangle)=1$
by (2). $\Box$  \\ \\ \\
{\bf Proof of Proposition \ref{prop1}}
Again we go through few steps. \\ \\
\underline{Step 1}. $[c,d,[a,b]]$ is in $Z(G)$. \\ \\
From Lemma \ref{lem1}, we know that $\langle a,a^{b},c,d\rangle $ is nilpotent of class at most $4$. Hence, using (2), we have 
$$\begin{array}{l}
\mbox{}[c,d,[a,b],c]=[c,d,a^{-1}a^{b},c]=[c,d,a^{-1},c][c,d,a^{b},c]=1. \\
\mbox{}[c,d,[a,b],d]=[c,d,a^{-1}a^{b},d]=[c,d,a^{-1},d][c,d,a^{b},d]=1.
\end{array}$$
By symmetry $[a,b,[c,d]]$, and hence $[c,d,[a,b]]$, commutes with $a$ and $b$. Therefore $[c,d,[a,b]]$ is in $Z(G)$. \\ \\
\underline{Step 2}. $[a,b,c,[d,a]]$ is in $Z(G)$. \\ \\ 
Using Lemma \ref{lem1} and (2), we have 
$$[a,b,c,[d,a]]=[a^{-1}a^{b},c,[d,a]]=[a^{-1},c,[d,a]][a^{b},c,[d,a]]=[a^{b},c,[d,a]].$$ 
By Lemma \ref{lem1} this commutes with $a,c$ and $d$. It remains to
see that it commutes with $b$. Now, using (1) and Step 1, $[c,b,a]^{[d,a]}=[[c,b][c,b,[d,a],a]=[c,b,a]$ and thus $[c,b,a]$ commutes with $[d,a]$. Hence, using again (1) and (2),
\begin{eqnarray*}
 \mbox{}[a,b,c,[d,a]] & = & [[c,b,a][a,c,b],[d,a]] \\ 
                                   & = & [a,c,b,[d,a]] \\
                                  & = & [a^{-1}a^{c},b,[d,a]] \\ 
                                & = & [a^{-1},b,[d,a]][a^{c},b,[d,a]] \mbox{\  \ (Using Lemma \ref{lem1})},
\end{eqnarray*}
that commutes with $b$ by Lemma \ref{lem1}. Hence $[a,b,c,[d,a]]$ is in $Z(G)$. \\ \\
\underline{Step 3}. $[a,b,c,d]$ is in $Z_{2}(G)$. \\ \\
We have
\begin{eqnarray*}
\mbox{}[a,b,c,d] & =& [a^{-1}a^{b},c,d]] \\
                             & = & [[a^{-1},c]^{a^{b}}[a^{b},c],d] \\
                             & = &  [[a^{-1},c][a^{b},c],d]\mbox{\ \ (as }[a^{-1},c]\mbox{ commutes with }a^{b}=a[a,b]\mbox{ by (1,2))} \\
                            & = & [a^{-1},c,d][a^{b},c,d] \mbox{\  \ (by Lemma \ref{lem1}).}
\end{eqnarray*}
By (2) this commutes with $c$ and $d$. Then we have modulo $Z(G)$, using (2) and Step 2,
     $$[a,b,c,d]^{a}=[a,b,c,d[d,a]]=[a,b,c,d]^{[d,a]}=[a,b,c,d]$$ 
and $[a,b,c,d]$ commutes with $a$. But, using (2), $[a,b,c,d]=[[b,a]^{-1},c,d]=[b,a,c]^{-1},d]=[b,a,c,d]^{-1}$, where for the last equality we used Lemma \ref{lem1}. By the previous calculations this element commutes with $b$. Hence $[a,b,c,d]$ is in $Z_{2}(G)$. \\ \\
\underline{Step 4}. $G$ is nilpotent of class at most $5$.  \\ \\
By Step 3, we know that modulo, $Z_{2}(G)$, $[a,b,c]$ commutes with $d$. By (2) it also commutes with $a,b$ and $c$. Hence, $[a,b,c]\in Z_{3}(G)$ and by symmetry $[a,b,d]\in Z_{3}(G)$. 
It follows that modulo $Z_{3}(G)$ we have that $[a,b]$ commutes with $c$ and $d$ and by (1) it also commutes with $a,b$. Hence $[a,b]$ is in $Z_{4}(G)$ and by symmetry $[a,c]$ and $[a,d]$ as well. Hence $a\in Z_{5}(G)$ and by symmetry $b,c$ and $d$ as well. It follows that $G$ is nilpotent of class at most $5$. $\Box$
\begin{prop} \label{prop2}
Let $X$ be a partial strong sandwich set in some group $G$ and let $a,b\in X$. Then $X\cup \{[a,b]\}$ is also a partial strong sandwich set in $G$. 
\end{prop}
{\bf Proof}\ \ Let $c,d\in X$ and $f,g\in \langle X\rangle$. It suffices to show that $\langle [a,b],c^{g}\rangle$ is nilpotent of class at most $2$ and that $\langle [a,b],c^{f},d^{g}\rangle$
is nilpotent of class at most $3$. Without loss of generality we can assume that $f=g=1$. As $\langle a,b,c\rangle$ is nilpotent of class at most $3$, the first assertion is immediate. By Proposition \ref{prop1} we know that $\langle a,b,c,d\rangle$ is nilpotent of class at most $5$. From this it follows that all commutators of weight $5$ in $[a,b],c,d$ with at least two occurrences of $[a,b]$ are trivial. It therefore only remains to see that $[[a,b],x_{1},x_{2},x_{3}]=1$ for all $x_{1},x_{2},x_{3}\in \{c,d\}$. But 
\begin{eqnarray*}
      [[a,b],x_{1},x_{2},x_{3}] & = & [a^{-1}a^{b},x_{1},x_{2},x_{3}] \\
                                               & = &  [a^{-1},x_{1},x_{2},x_{3}][a^{b},x_{1},x_{2},x_{3}] \mbox{\ \ (using Lemma \ref{lem1})} \\
                                              & = & 1 \mbox{\ \ (by (2)).} 
\end{eqnarray*}
This finishes the proof. $\Box$
\begin{prop}\label{prop3}
Let $X$ be a strong sandwich set in some group $G$ and let $a,b\in X$. Then $X\cup \{[a,b]\}$ is also a strong sandwich set. 
\end{prop}
{\bf Proof}\ \ Let $c,d\in X$ an $f,g\in \langle X\rangle $. Firstly it is immediate as before that $\langle [a,b],c^{g}\rangle $ is nilpotent of class at most $2$. That the same is true for 
$\langle [a,b],[a,b]^{g}\rangle$ follows from the fact that $\{ a,b,a^{g},b^{g}\}$ is a partial strong sandwich set and thus $\langle a,b,a^{g},b^{g}\rangle$ nilpotent of class at most $5$ by Proposition \ref{prop1}. It remains to show that $\langle [a,b],c^{f},d^{g}\rangle$, $\langle [a,b],[a,b]^{f},c^{g}\rangle$ and $\langle [a,b], [a,b]^{f}, [a,b]^{g}\rangle$ are nilpotent of class at most $3$. By Proposition \ref{prop2}, we know that $\{[a,b],c^{f},d^{g}\}$ is a partial sandwich set and thus $\langle [a,b],c^{f},d^{g}\rangle$ nilpotent of class at most $3$. \\ \\
Next consider the set $\{a,b,a^{f},b^{f}, c^{g}\rangle$. As this is a partial strong sandwich set, if follows by two applications of Proposition \ref{prop2} that $\{[a,b],[a,b]^{f},c^{g}\rangle$
is a partial strong sandwich set and therefore nilpotent of class at most 3. Finally consider the partial strong sandwich set $\{a,b,a^{f},b^{f},a^{g},b^{g}\}$. By three applications of Proposition \ref{prop2} we see that $\{[a,b],[a,b]^{f},[a,b]^{g}\}$ is also a partial strong sandwich set and therefore $\langle [a,b], [a,b]^{f}, [a,b]^{g}\rangle $ nilpotent of class at most $3$. 
$\Box$
\section{Nilpotence of finitely generated strong sandwich groups}
The rest of the proof of Theorem 1.1 is very similar to the corresponding part of the proof of the main result in \cite{Jab}. We include the proof for the convenience of the reader and also because the setting is a bit different. \\ \\
Let $X=\{x_{1}, x_{2},\ldots ,x_{r}\}$ be a strong sandwich set. In this section we prove that $G=\langle X\rangle$ is nilpotent. 
Let $\overline{X}$ be the closure of $X$ with respect to the commutator action. In other words $\overline{X}$ consists of all commutators in $X$ (in any order and with any bracketing). It follows by iterated use of Proposition \ref{prop3} that $\overline{X}$ is a strong sandwich set. 
\begin{lemm}
Let $u,v,w\in \overline{X}$. Then $[u,[v,w]]=[u,v,w][u,w,v]^{-1}$. 
\end{lemm}
{\bf Proof }\ \ As $\langle u,v,w\rangle $ is a strong sandwich group, we know that it is nilpotent of class at most $3$. The result now follows from this and the Hall-Witt identity. \\ \\
Our proof makes use of the notion of standard words (see for example \cite{Vaug}) which played a crucial role in Chanyshev's proof of the theorem on sandwich algebras by Kostrikin and Zel'manov \cite{Kos}. Our proof resembles the work of Chanyshev in outline although with group commutators instead of works in a Lie ring and it is similar to the proof in \cite{Jab}. \\ \\
{\bf Standard words.} Let $x_{1},\ldots ,x_{r}$ be free variables. Consider the set $A$ of all words
                                            $$x_{i(1)}\cdots x_{i(n)},\ 1\leq i(1),\ldots ,i(n)\leq r,\ n\geq 0.$$
We order these words as follows: $x_{i(1)}\cdots x_{i(n)}<x_{j(1)}\cdots x_{j(m)}$ if either for some $t<\mbox{min}\{m,n\}$ we  have $x_{i(1)}=x_{j(1)}, \ldots ,x_{i(t)}=x_{j(t)}$ and $x_{i(t+1)}<x_{j(t+1)}$, or $m<n$ and $x_{j(1)}=x_{i(1)},\ldots ,x_{j(m)}=x_{i(m)}$. This gives us a total order on $A$. \\ \\
{\bf Definition}. We say that a word $x_{i(1)}\cdots x_{i(n)}$ is standard if for all $2\leq t\leq n$ we have $x_{i(t)}\cdots x_{i(n)}x_{i(1)}\cdots x_{i(t-1)}<x_{i(1)}\cdots x_{i(n)}$. \\ \\
We make use of the following property for standard words. If $c$ is a standard work of length at least $2$, then $c=ab$ for some standard words $a$,$b$ where $a>b$. Among all such decompositions we pick the one where $a$ is the largest. Although the choice of $a$ is irrelevant in what follows. \\ \\
{\bf Definition} To each standard word $c$ we associate a group commutator, denoted $[c]$, recursively as follows. Firstly $[x_{i}]=x_{i}$ for $i=1,\ldots ,r$. Then if the length $l(c)$ of $c$ is at least $2$ and $c=ab$ for standard words $a,b$, then $[c]=[[a],[b]]$. \\ \\
{\bf Definition}. To each word $c=x_{i(1)}\cdots x_{i(n)}$ in $A$ we associate the left normed commutator $\mbox{com}(c)=[x_{i(1)},\ldots ,x_{i(n)}]$. \\ \\
{\bf Theorem 1.1} {\it $G=\langle X\rangle$ is nilpotent}. \vspace{2mm} \\
{\bf Proof}\ \ 
Let 
                     $$H=Z_{\infty}(G)=\bigcup_{i=0}^{\infty}Z_{i}(G)$$
be the hyper centre of $G$. To show that $G$ is nilpotent it suffices to show that $\gamma_{m}(G)\leq H$ for some positive integer $m$. We argue by contradiction and suppose this is not the case. We then get an infinite sequence $(x_{\alpha(i)})_{i=1}^{\infty}$ where $u_{n}=x_{\alpha (1)}\cdots x_{\alpha(n)}$ is the smallest word in $A$ of length $n$ such that $\mbox{com}(u_{n})\not\in H$. Before proceeding we need a lemma.
\begin{lemm} 
Let $n\geq 1$ and $c$ a standard word of length $m$. \vspace{2mm}\\
(1) If $u_{n}c<u_{n+m}$, then $[\mbox{com}(u_{n}),[c]]\in H$. \\
(2) If $u_{n}c=u_{n+m}$, then $[\mbox{com}(u_{n}),[c]]H=\mbox{com}(u_{n+m})H$. 
\end{lemm}
{\bf Proof}\ \ We prove this by induction on $m$. For $m=1$, the statement (2) is obvious while (1) follows directly from our choice of the sequence $(x_{\alpha(i)})_{i=1}^{\infty}$. Let $m\geq 2$ and suppose (1) and (2) hold for smaller values of $m$. Consider first (1). Let $c=ab$, where $a,b$ are standard words of lengths $s,t$. Then 
 \begin{eqnarray*}
                [\mbox{com}(u_{n}),[c]] & = & [\mbox{com}(u_{n}),[[a],[b]]] \\ 
                                                         & \stackrel{L3.1}{=} & [\mbox{com}(u_{n}),[a],[b]]\cdot [\mbox{com}(u_{n}),[b],[a]]^{-1}.
\end{eqnarray*}
As $u_{n}ab=u_{n}c<u_{n+m}$, we must have $u_{n}a\leq u_{n+s}$. If $u_{n}a<u_{n+s}$, then we have $[\mbox{com}(u_{n}),[a]]H=H$ by the induction hypothesis. If on the other hand
$u_{n}a=u_{n+s}$, then $u_{n+s}b<u_{n+m}$ and thus by the induction hypothesis 
           $$[\mbox{com}(u_{n}),[a],[b]]H=[\mbox{com}(u_{n+s}),[b]]H=H.$$
As $u_{n}ab=u_{n}c<u_{n+m}$, we must have $u_{n}a<u_{n+s}$. If $u_{n}a<u_{n+s}$, then we have $[\mbox{com}(u_{n}),[a]]H=H$ by the induction hypothesis. If on the other hand $u_{n}a=u_{n+s}$ then $u_{n+s}b<u_{n+m}$ and thus by the induction hypothesis 
                                   $$[\mbox{com}(u_{n}),[a],[b]]H=[\mbox{com}(u_{n+s}),[b]]H=H.$$
As $c=ab$ is standard, we have $ba<c\leq u_{n+m}$. The same argument as before gives $[\mbox{com}(u_{n}),[b],[a]]H=H$. Hence it follows that $[\mbox{com}(u_{n}),[c]]H=H$ and we have proved (1) for $m$. \\ \\
We next turn to (2). As $u_{n}ba<u_{n}ab=u_{n}c=u_{n+m}$, the same argument as above shows that $[\mbox{com},[b],[a]]H=H$. As $u_{n}a=u_{n+s}$ and $u_{n+s}b=u_{n+m}$, it thus follows from the induction hypothesis that 
\begin{eqnarray*}
         [\mbox{com}(u_{n}),[c]]H & = & [\mbox{com}(u_{n}),[a],[b]]H \\
                                                & = & [\mbox{com}(u_{n+s}),[b]]H \\
                                                & = & [\mbox{com}(u_{n+m})H.
\end{eqnarray*}
This finishes the proof of lemma. $\Box$ \\ \\
{\bf Continuation of the proof of Theorem 1.1}. By Theorem 3.1.10 in \cite{Vaug} we know that there exists and integer $N(r)$ such that any word in $A$ of length greater than $N(r)$ must contain a subword $cc$, $cx_{i}c$ or $x_{i}cx_{i}$ where $c$ is a standard word. In particular, the word $x_{\alpha(2)}\cdots x_{\alpha(2+N(r))}$ must contain one of these. If $c$ has length $m$, one of the following mus therefore hold 
  $$u_{n+2m}=u_{n}cc;\ u_{n+2m+1}=u_{n}cx_{i}c;\mbox{ or }u_{n+m+2}=u_{n}x_{i}cx_{i}$$
for some positive integers $n,m$. Now using the fact that $\langle \mbox{com}(u_{n}), [c],x_{i}\rangle$ is a strong sandwich group, and therefore nilpotent of class at most $3$, we see in the first case, using Lemma 3.2, that 
                     $$H=[\mbox{com}(u_{n}),[c],[c]]H=[\mbox{com}(u_{n+m}),[c]]H=\mbox{com}(u_{n+2m})H.$$
This gives the contradiction that $\mbox{com}(u_{n+2m})in H$. Similarly for the other cases we get contradictions from 
$$\begin{array}{l}
          H=[\mbox{com}(u_{n}),[c],x_{i},[c]]H=\mbox{com}(u_{n+2m+1})H, \\
         H=[\mbox{com}(u_{n}),x_{i},[c],x_{i}]H=\mbox{com}(u_{n+m+2})H.
\end{array}$$
From these contradictions we conclude that $G$ must be nilpotent. This finishes the proof of Theorem 1.1. $\Box$ \\ \\
As we pointed out in the introduction, Theorem 1.2 follows from this. \\ \\
{\bf Remark}. It follows from this that when all groups satisfying a law $w=1$, in some variables $x_{1},\ldots ,x_{n}$, are locally nilpotent, then the new variety, satisfying $[x_{n+1},w,w,w]=1$ and $\langle w,w^{x_{n+2}},w^{x_{n+3}}\rangle$  being nilpotent of class at most $3$, is also locally nilpotent. Starting for example with $w=x_{1}^{4}$, this gives us a sequence of locally nilpotent varieties.

\end{document}